\documentclass[conference,10pt]{IEEEtran}
\usepackage{listings}
\usepackage{balance}
\usepackage{graphicx}
\usepackage{eqnarray,amsmath,mathtools, nccmath}
\usepackage{amsfonts}
\newcommand{\Z}{\mathbb{Z}}

\DeclarePairedDelimiter{\nint}\lfloor\rceil

\IEEEoverridecommandlockouts
\allowdisplaybreaks[4]

\begin{document}

\title{\Large \bf
OPTIMAL ROTATIONAL LOAD SHEDDING \\ VIA BILINEAR INTEGER PROGRAMMING
}

\author{Atif Maqsood, Yu Zhang, and Keith Corzine \\
Department of Electrical and Computer Engineering \\
University of California, Santa Cruz
\thanks{Emails: {\tt\{amaqsood,zhangy,corzine\}@ucsc.edu.}}
\thanks{This work was supported by the Faculty Research Grant (FRG) of UC Santa Cruz.}
}

\maketitle

\begin{abstract}
This paper addresses the problem of managing rotational load shedding schedules for a power distribution network with multiple load zones. An integer optimization problem is formulated to find the optimal number and duration of planned power outages. Various types of damage costs are proposed to capture the heterogeneous load shedding preferences of different zones. The McCormick relaxation along with an effective procedure feasibility recovery is developed to solve the resulting bilinear integer program, which yields a high-quality suboptimal solution. Extensive simulation results corroborate the merit of the proposed approach, which has a substantial edge over existing load shedding schemes.
\end{abstract}

\section{Introduction}
Over the last few decades, power consumption has been steadily increasing around the world \cite{power}. This can be attributed to a rise in availability of electrical and electronic devices as well as social and cultural factors that promote our dependence on those machines. For many developing countries, the fast population growth is the main driver of increased energy consumption  \cite{popu}. In addition, ludicrous electricity consumption has been recently observed due to large numbers of geo-distributed data centers and bitcoin mining machines \cite{bitcoin}. 

The increased demand poses a heavy burden on our aging and stressing energy infrastructure. Consider Pakistan as an example, which has the world's sixth largest population (\mbox{197 million}), and two of the top twelve most populated cities in the world \cite{pak}. With such a high population and limited resources, it is no surprise that the demand for power in recent years has consistently exceeded the generation capacity leading to shortfalls and blackouts \cite{pakfall,impact}. Efforts to increase the power generation and replace the aging infrastructure to support higher levels of demand have been slow and hampered by various financial and political factors \cite{USIP}.

As the utility company in Lahore, Pakistan, Lahore electric supply company (LESCO) has been conducting the so-called rotational load shedding in the city to deal with frequent shortfalls~\cite{lesco}. Load shedding is a premeditated power outage in some load zones when it is expected that the power supply cannot meet the demand. Rotational load shedding ensures that over the time horizon planned power outages will be shared across different zones in a pre-determined fashion. Refusing to implement rotational load shedding can result in unscheduled outages bringing major inconvenience, and even catastrophic consequences like large-scale rolling blackouts. Infrastructure such as distribution lines and transformers may endure permanent damages due to overheating. It is therefore preferred to shed load rather than take a risk of undesirable outcomes \cite{collapse}.

\begin{figure}[!t]
\centering
\includegraphics[width=3.4in]{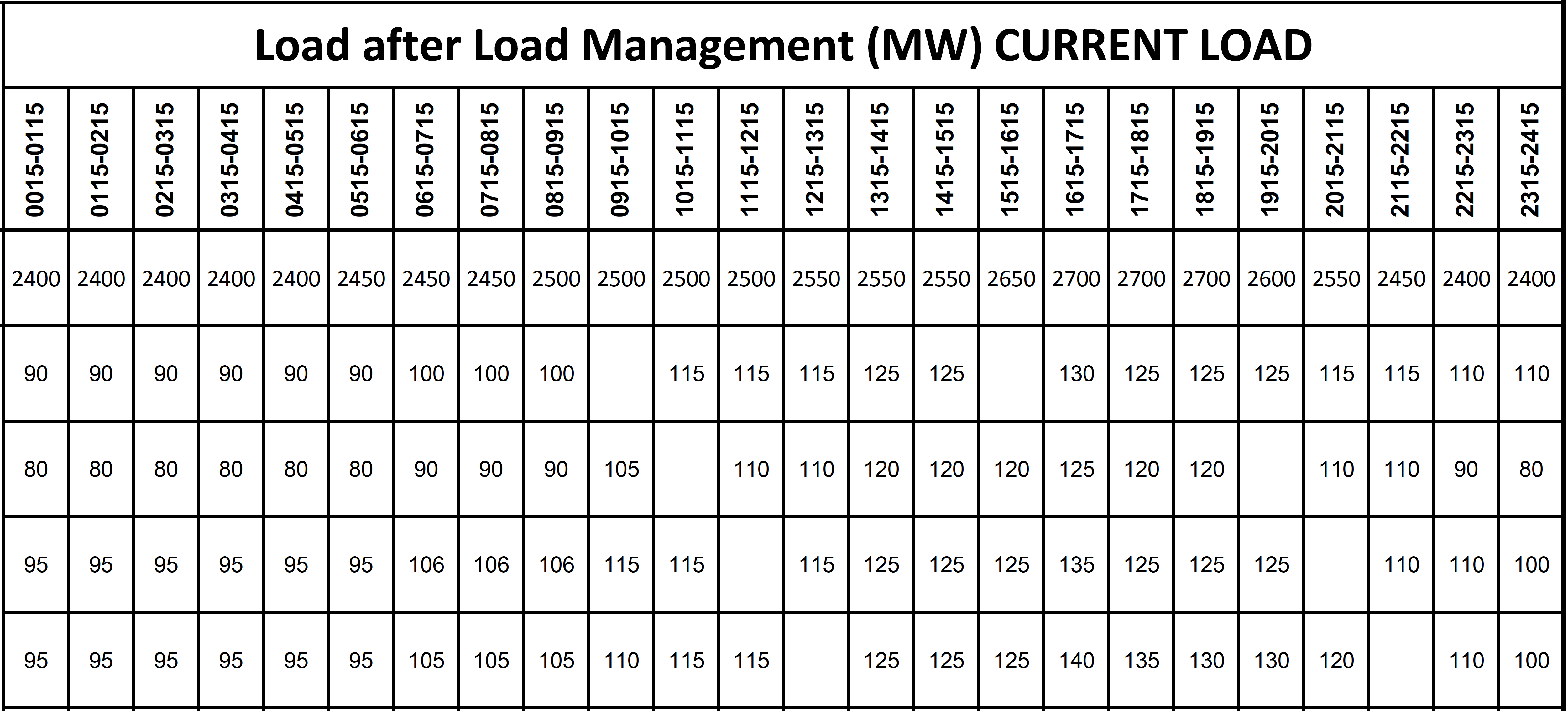}
\caption{Load management program for Lahore, Pakistani: As per AT\&C losses W.E.F 01-05-2018 \cite{lesco}.}
\label{lesco1}
\vspace*{-0.3cm}
\end{figure}

\begin{figure}[!t]
\centering
\includegraphics[width=3.4in]{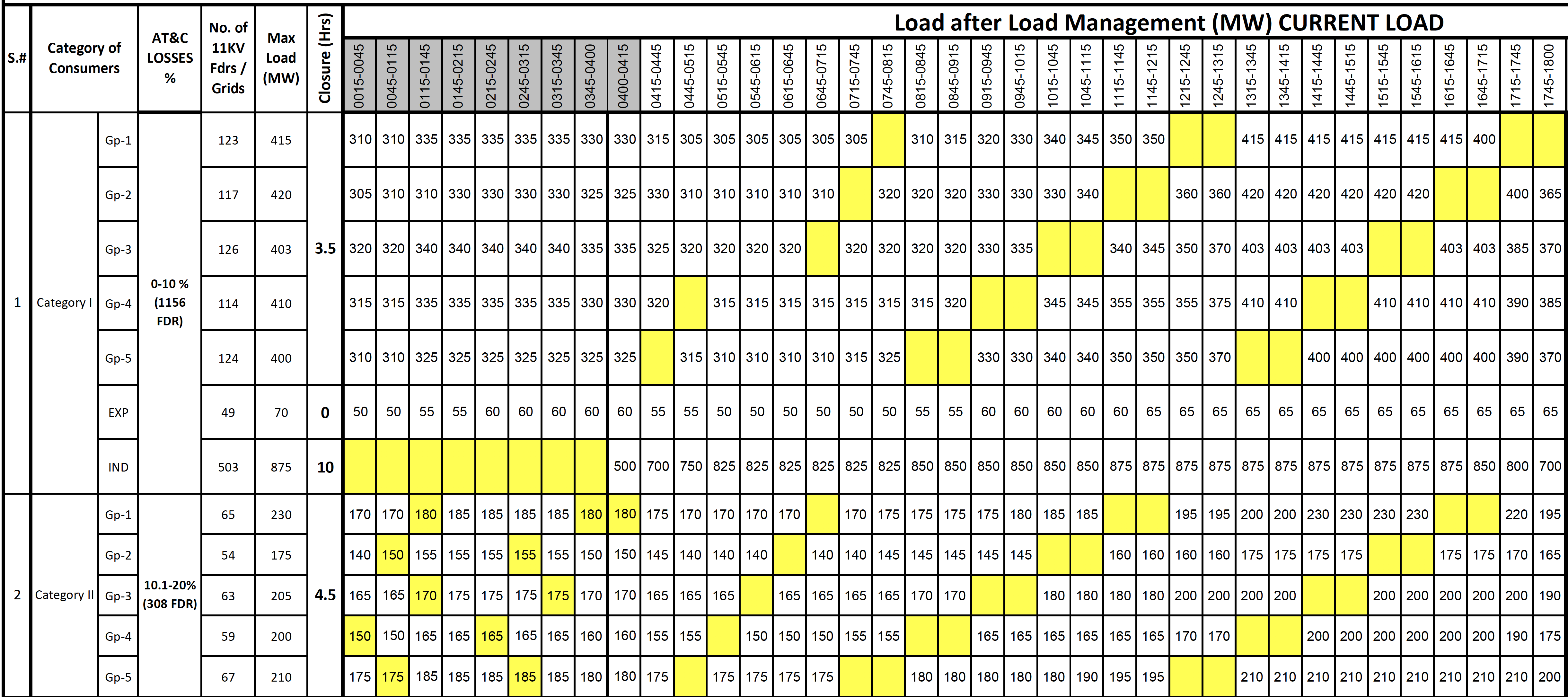}
\caption{Lesco load management program for Lahore Ramadan 2018 on the basis of AT\&C losses with no load management during Seher, Iftar and Travih timings \cite{lesco}.}
\label{lesco2}
\vspace*{-0.3cm}
\end{figure}

According to LESCO's most recent announcement, there is a load shedding of 12 hours per day that is currently implemented \cite{lesconews}. Over the last decade this number has been varying from 10 to 18 hours per day. Figure \ref{lesco1} shows the real data for the load shedding schedules in Jan 2018 \cite{lesco}. Over 1500 feeders in the city are classified into five categories, some of which are further divided into different groups. Each group experiences a different schedule and total outage hours. The duration of load shedding varies from 2 to 6 hours per day. At the moment, there are no standards for such a division and schedule variations for different zones. Hence, it is completely up to the utility and local government to decide how to allocate the unavoidable shortfalls among its customers. Figure \ref{lesco2} shows the load shedding schedules for the lunar month of Ramadan in 2018, where the outage patterns also varies based on the type of load; i.e., industrial vis-\`{a}-vis residential load zones. Note that even for residential zones, the load shedding varies from 3.5 to 7 hours per day due to the increased power consumption in summer.

Many developing countries have also been dealing with rolling blackouts via rotational load shedding. Egypt has a shortfall every summer since 2010 because of the high power demand and increased use of air conditioners and fans \cite{egypt}. Ghana has had year-long shortfalls and load shedding since 2012, sometimes up to 12 hours a day \cite{ghana}; see also the power crisis in India and South Africa \cite{rsa,indiafall}. The load shedding hours recorded over a period of 3 months in Chennai are given in \cite{chennai}. The data show that  there is no consistency of outage but obvious pattern in the start time of power interruption that indicates peak hours. Also the outage time is inconsistent. In the United States, e.g., Northern California, the investor-owned utility Pacific Gas and Electric Company (PG\&E) can conduct rotating outages if the ISO deems it necessary \cite{pge}. High power demand is not the only reason for shortfall. Unexpected disasters that lead to failure in power generation sometimes can result in rolling blackouts and rotational load shedding; e.g., the 2011 Tohoku tsunami~\cite{japan}, the 2014 Ukraine crisis, and the 2015 North American heat wave. Clearly, the long-term solution for energy shortfall is to increase the generation capacity, upgrade existing power infrastructure, and incentivize customers to reduce power usage during peak hours via demand response~\cite{dr}. However, rotational load shedding will still be a viable approach for some countries and regions to cope with shortfalls in the next few decades.


Focusing on this problem, the present paper proposes an optimization framework aiming to minimize the total damage cost while respecting various load shedding preferences among different zones (see Section~\ref{sec:formulation}). A relaxation technique is leveraged to linearize the bilinear terms, and yield a high-quality suboptimal solution (Section~\ref{sec:solve}).  
Section~\ref{sec:test} reports the numerical performance of the novel approach. Finally, concluding remarks are given in Section~\ref{sec:summary}.

\section{Problem Formulation} \label{sec:formulation}       
Consider a power distribution system consisting of $N$ load zones, which are indexed by $n \in \mathcal{N} :=\{1,2,\ldots, N\}$. We assume that each load zone can be shed independently without affecting the rest of the system.  
Let $k_n$ and $d_n$ denote the number of outages, and the duration (in hours) of each outage for the $n$-th zone, respectively. In addition, each zone has a  cost function $C_n(k_n, d_n)$, which essentially reflects the damage cost due to the load shedding.

\subsection{Heterogeneous Load Shedding Costs}
The modeling of the load shedding damage cost functions requires the knowledge of end users' consumption habits, as well as  financial ramifications due to the loss of load. An accurate modeling can be challenging since the inconvenience caused by load shedding sometimes is not necessarily related to financial losses. This however is out of the focus of the present work. 

In this paper, we first classify all load zones as residential, commercial, and industrial areas that feature heterogeneous load shedding costs.
Industrial zones often prefer few outages of longer duration, rather than frequent outages of short duration. This is due to the high cost of shut down and restart operations, especially when large machines are involved. In contrast,  residential zones would prefer the other way because most residence and offices rely on backup power such as UPS, which can last for a few hours. Such practical preferences are corroborated by the load shedding program in Figure \ref{lesco2}; see also the planned rolling blackouts in South Africa 2015 for similar preferences \cite{rsa2}. 
Specifically, the damage costs are modeled as increasing functions in $k_n$ and $d_n$, which may serve as approximation of the true costs among zones. 
For example, an industrial zone 1 and a residential zone 2 can take the following damage costs respectively
\begin{subequations}
\begin{align}
C_1(k_1,d_1) &=1000k_1+3000d_1k_1,\\
C_2(k_2,d_2) &=1000d_2+100d_2k_2\, .
\end{align}
\end{subequations}
The functions defined above show the relative preference that each zone places on their load shedding schedules. Zone 1 assigns a high weight cost for avoiding a large number of frequent outages while zone 2 gets high cost for outages with long duration. Note that both cost functions have the product term 
$d_n \times k_n$ which represents the total amount of time for loss of load. In this example zone 1 has a higher coefficient penalizing the unavailability of power. This may reflect the fact that some utility companies would prefer to shed power in residential zones rather than industrial ones. 

A more general cost function is given as 
\begin{align}
C_n(k_n,d_n) := a_{n,1}d_nk_n+a_{n,2}d_n+a_{n,3}k_n, 
\end{align}
where the coefficients  $a_{n,1}, a_{n,2}, a_{n,3} \in \mathbb{R}_+$ reflect the load shedding preference of each zone.  
These cost functions can also vary with time of the day, season, and year. An outage of one hour may have different ramifications in summer and winter, as well as during the day and night. Here we simply consider the worst case damage for the period considered.

The objective of the utility is to find the optimal $d_n$ and $k_n$ for each zone such that the overall damage cost can be minimized.
We can simply consider minimizing the cumulative cost for all zones as given by 
\begin{align}
C_{\mathrm{tot}}=\sum_{n=1}^{N} C_n. \label{costtotal}
\end{align} 
More complicated composite cost functions are applicable if the load shedding of a zone intensifies the damage in its neighborhood zones. 

\subsection{Operational Constraints}
Without any operational constraints the optimal decision would be to have no outages for any zones. This implies that we need to supply more power demand than the total generation capacity.  The long lasting overburdened infrastructure may eventually lead to deterioration of voltage or frequency and rolling blackouts.  To avoid such dangerous consequences, we need to ensure that energy being shed is no less than the expected short fall, which is denoted by $E_{\mathsf{sf}}$. This hard constraint can then be represented as
\begin{align}
\sum_{n=1}^N d_nk_nP^{\mathrm{ave}}_n  \geq E_{\mathsf{sf}} \label{loadshed}
\end{align}
where $P^{\mathrm{ave}}_n$ is the estimate of the average power consumed in zone $n$. 
Note that the expected short fall is usually independent of the desired reliability, which can be obtained by predicting the power generation and demand over a period of time; e.g. a season or a year. Estimate of the average power $P^{\mathrm{ave}}_n$ can be obtained by widely used load forecasting methods; see \cite{pred}. Apart from this operational condition, due to the practical limits of outage duration $d_n$ and  number of outages $k_n$, we also have the following two box constraints
\begin{align}
0 &\leq k_n \leq \overline{k}_n,\quad \forall n \in \mathcal{N} \\
\underline{d}_n &\leq d_n \leq \overline{d}_n, \quad \forall n \in \mathcal{N}.
\end{align} 

To promote the fairness among all zones, we further set a limit $\overline{C}_\delta$ for the cost difference between each pair of two adjacent zones, as given by the following constraint
\begin{align}
| &C_n(k_n,d_n)-C_{n+1}(k_{n+1},d_{n+1})| \leq \overline{C}_\delta,\quad \forall n \in \mathcal{N}. \label{costdiff}
\end{align}

\section{Solving the Bilinear Integer Program}\label{sec:solve}    
The minimum duration of each power outage is set to be 15 minutes. Hence, the duration of each planned load shedding becomes $0.25d_n$ while $d_n$ is an integer number. To this end, the task of optimal rotational load shedding can be formulated as an optimization problem that minimizes the total load shedding cost \eqref{costtotal} subject to all operational constraints \eqref{loadshed}-\eqref{costdiff}. The difficulty in solving the resulting optimization problem is two-fold: i) this is an integer program since both decision variables $\{d_n, k_n\}_{n\in \mathcal{N}}$ are integer numbers; and ii) bilinear terms $\{d_n\times k_n\}_{n\in \mathcal{N}}$ are involved in the objective and some constraints. Hence, off-the-shelf integer programming solvers are not directly applicable. 

\subsection{McCormick relaxations}
To deal with such a challenge, we leverage the McCormick relaxations, and linearize each bilinear term by introducing a new variable $w_n := d_n\times k_n$; see e.g., \cite{milp}. Interestingly, for our particular problem, the new variable $w_n$ represents the total amount of time for planned load shedding in zone $n$; i.e., unavailability of power in hours. Replacing all bilinear terms with the introduced new variables, the relaxed problem is given as 
\begin{subequations}\label{relaxP}
\begin{align}
\hspace{-0.6cm}\min\limits_{\{w_n,d_n,k_n\}}\,\,  &\sum_{n=1}^N a_{n,1}w_n+\frac{1}{4}a_{n,2}d_n+a_{n,3}k_n \label{obj} \\
\mathrm{s.t.} \quad\,\, &\frac{1}{4}\sum_{n=1}^N w_nP^{\mathrm{ave}}_n-E_{\mathsf{sf}} \geq 0 \\
& 0 \leq k_n \leq \overline{k}_n , \, \forall n \in \mathcal{N} \label{krange} \\
&\underline{d}_n \leq d_n \leq \overline{d}_n , \, \forall n \in \mathcal{N} \\
&\underline{d}_nk_n \leq w_n \leq \overline{d}_nk_n , \, \forall n \in \mathcal{N}  \label{wrange1}\\
&\overline{d}_nk_n + \overline{k}_nd_n-\overline{k}_n\overline{d}_n \leq w_n \leq \overline{k}_nd_n , \, \forall n \in \mathcal{N} \label{wrange2}  \\
& \bigl |a_{n,1}w_n+\frac{1}{4}a_{n,2}d_n+a_{n,3}k_n -a_{n+1,1}w_{n+1} \notag \\
&+\frac{1}{4}a_{n+1,2}d_{n+1}+a_{n+1,3}k_{n+1}\bigr | \leq \overline{C}_\delta, \, \forall n \in \mathcal{N} \label{cons} \\
& w_n,d_n,k_n \in \Z_+, \, \forall n  \in \mathcal{N}\,.
\end{align}
\end{subequations}
Note that the lower and upper limits of $w_n$ [cf. \eqref{wrange1}-\eqref{wrange2}] are derived by the limits of $d_n$ and $k_n$ via the McCormick envelope. 
The factor $1/4$ is introduced due to the fact that we set 15 minutes as the base unit of the outage time, and each outage duration is $d_n$ times of the base.

\subsection{Feasible Solution Recovery}
By using the McCormick relaxation, problem \eqref{relaxP} becomes an integer program with linear objective function and constraints, which can be dealt with integer solvers such as {\tt{Cplex}} and {\tt{Gurobi}}. Let $\{k_n^*,d_n^*,w_n^*\}_{n\in \mathcal{N}}$ denote an optimal solution to the relaxed problem \eqref{relaxP}. If $w_n^* = d_n^* k_n^*$ for $n = 1,2,\ldots, N$, then the relaxation is exact, and we essentially solve the original load shedding problem with the bilinear terms. If the relaxation is inexact, the optimal solution $d_n,k_n$ violates the minimum short fall constraint, i.e.
$
\sum_{n=1}^N d_n^*k_n^* P^{\mathrm{ave}}_n  < E_{\mathsf{sf}}.
$

In this case, we need an extra procedure to recover a feasible solution to the original problem. We multiple all outages with a factor such that the relative number of outages between different zones stay the same while the overall number of outages is increased to meet the shortfall. That is, i) keeping $d_n^*$ unchanged; and ii) rescaling $k_n^*$ and rounding it to the nearest integer
\begin{align}
\check{d}_n^* &= d_n^* \\
\check{k}_n^* &= \nint*{\frac{4E_{\mathsf{sf}}}{\sum_{n=1}^N d_n^*k_n^*P^{\mathrm{ave}}_n}k_n^*}\,.
\end{align}

\section{Numerical Tests}\label{sec:test}  
\begin{figure}[!t]
\centering
\includegraphics[width=3.4in]{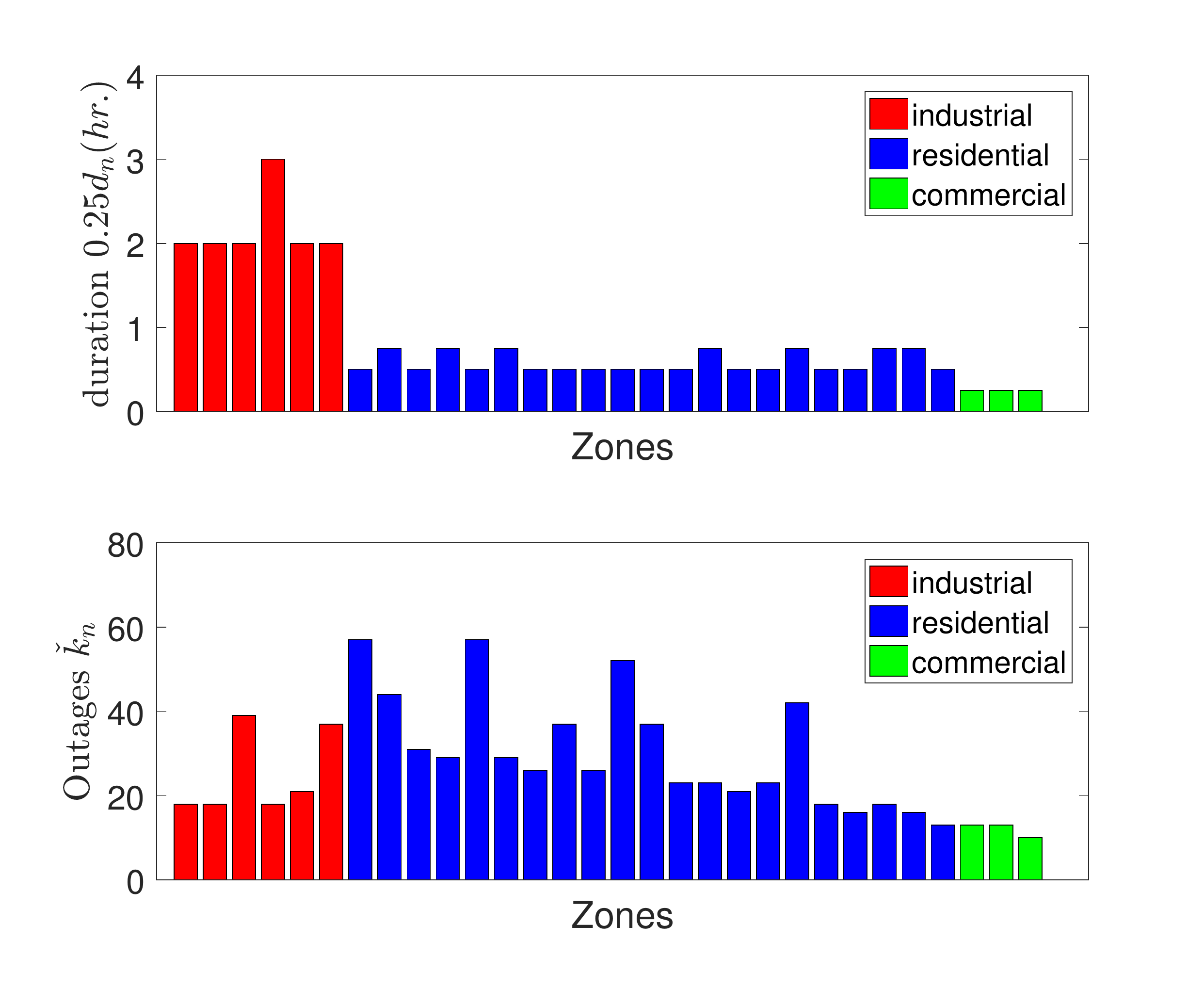}
\caption{Optimal load shedding schedules.}
\label{result}
\vspace*{-0.4cm}
\end{figure}

The resulting problem \eqref{relaxP} is solved by the {\tt{Gurobi}} along with {\tt{CVX}}. The simulation parameters are summarized in Table \ref{para} and \ref{para2}. The total period considered in the simulation is 30 days. We have three types of load zones, each of which has a unique cost function and average power $P^{\mathrm{ave}}_n$. This is a slightly more complicated system than the LESCO example discussed in the introduction. LESCO provides power to 1747 zones that are  grouped into 17 load shedding schedules (cf. Figure \ref{lesco2}). Figure \ref{result} are the results of optimal load shedding results. The top plot shows the duration $d_n$ per outage, which are multiples of the base unit time 15 minutes. The variance among different categories reflect the preferences given in Table \ref{para2}. It can be seen that the first 6 zones of industrial areas have generally a longer duration than other zones. The optimal   numbers of outages $\check{k}_n$ shown in the bottom plot are within the constraints for all zones. Due to the associated high costs, three commercial zones get the smaller numbers of outages compared with the other types of zones. 
\begin{table}[]
\centering
\caption{Simulation parameters for CVX.}
\label{para}
\begin{tabular}{|c||c|c|}
\hline
Parameter & Description & Value \\ \hline \hline
$n$ & Number of total zones & 30 \\ \hline
$n_{ind}$ & Number of industrial zones & 6 \\ \hline
$n_{res}$ & Number of residential zones & 21 \\ \hline
$n_{com}$ & Number of commercial zones & 3\\ \hline
 & Daily average power/hr &  \\ 
$P^{\mathrm{ave}}_n$ & for zone n & $500$ -- $1000$ MW\\ \hline
 & Deficit between Supply and & \\
$E_{\mathsf{sf}}$ & Demand over 30 days period & $5 \times 10^5$ MWh \\ \hline
 & Limit for cost differences & \\
$\overline{C}_\delta$ & among adjacent zones & 500 units \\ 
\hline
\end{tabular}
\end{table}
\begin{table}[]
\centering
\caption{Simulation parameters for load shedding cost functions.}
\label{para2}
\begin{tabular}{|c||c|c|c|}
\hline 
Parameter & Industrial & Residential & Commercial \\ \hline
\hline 
$a_{n,1}  $ & 50 - 150 & 50 - 150 & 500 - 600 \\ \hline 
$a_{n,2}  $ & 20 - 70 & 70 - 120 & 70 - 120 \\ \hline
$a_{n,3}  $ & 70 - 120 & 20 - 70 & 70 - 120 \\ \hline
$\overline{k}_n $ (per month)  & 50 & 200 & 20 \\ \hline
$\overline{d}_n  $ (hour per outage) & 4 & 2 & 0.5 \\ \hline
$\underline{d}_n $ (hour per outage) & 2 & 0.5 & 0.25 \\ \hline
\end{tabular}
\end{table}

\begin{figure}[!t]
\centering
\includegraphics[width=3.4in]{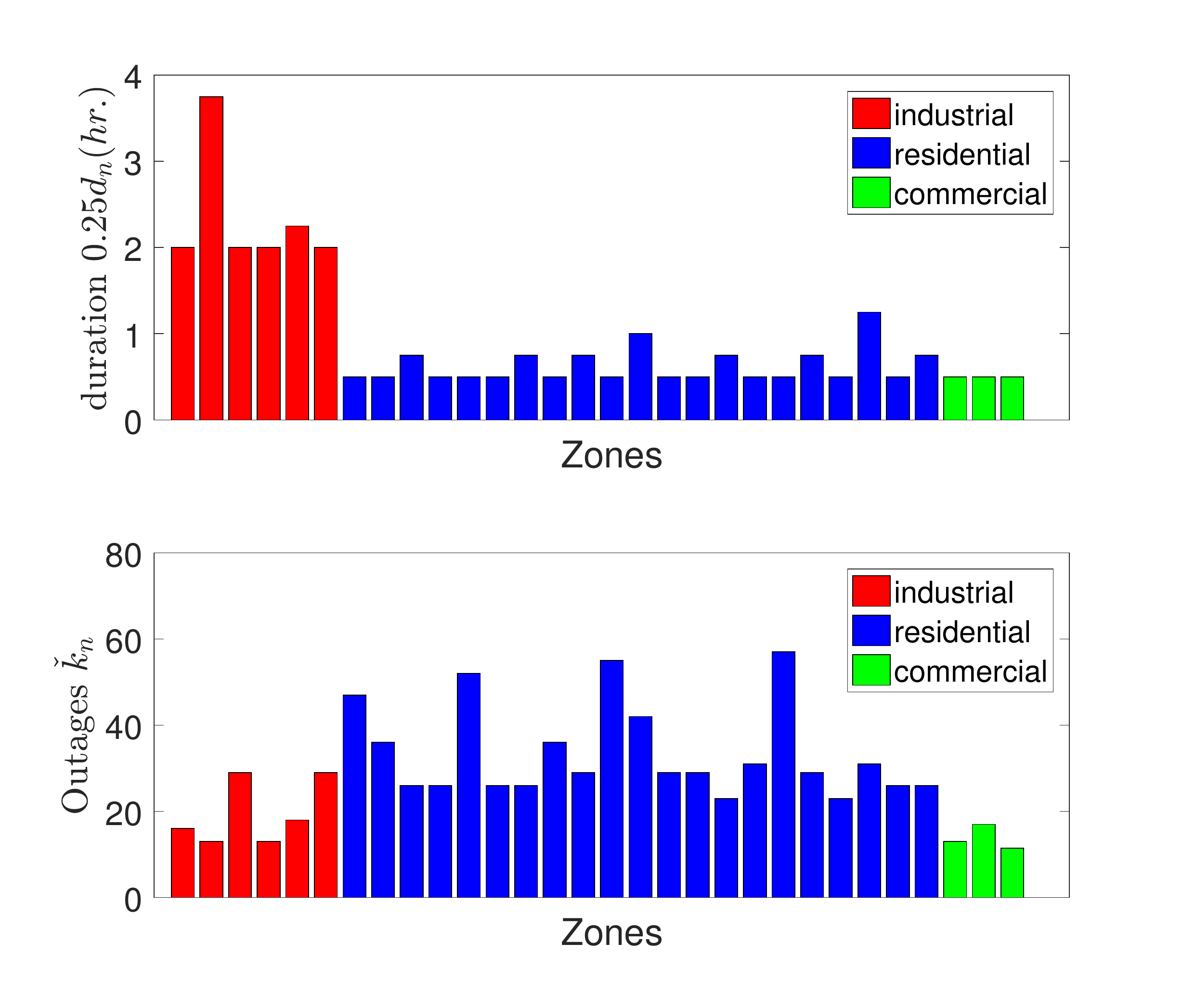}
\caption{Optimal load shedding schedules with a smaller $\overline{C}_\delta$.}
\label{result2}
\vspace*{-0.4cm}
\end{figure}

Thanks to the limited cost difference constraints \eqref{costdiff}, residential zones that are close to commercial zones have much less planned outages than those near industrial zones, as shown in Figure \ref{result}.
Figure \ref{result2} shows the optimal load shedding schedules where $\overline{C}_\delta = 100$ while all other parameters are kept the same as before. Clearly, the costs among different zones become more uniformly shared, compared with the results in \ref{result}. A smaller value of the limit $\overline{C}_\delta$ tightens the feasible set of \eqref{relaxP}, which yields a higher optimal cost.  Note that there exists a lower limit on $\overline{C}_\delta$ beyond which the optimization problem may become infeasible. The original value $\overline{C}_\delta = 500$ is used for the rest of the simulations.
  
\begin{figure}[!t]
\centering
\includegraphics[width=3.4in]{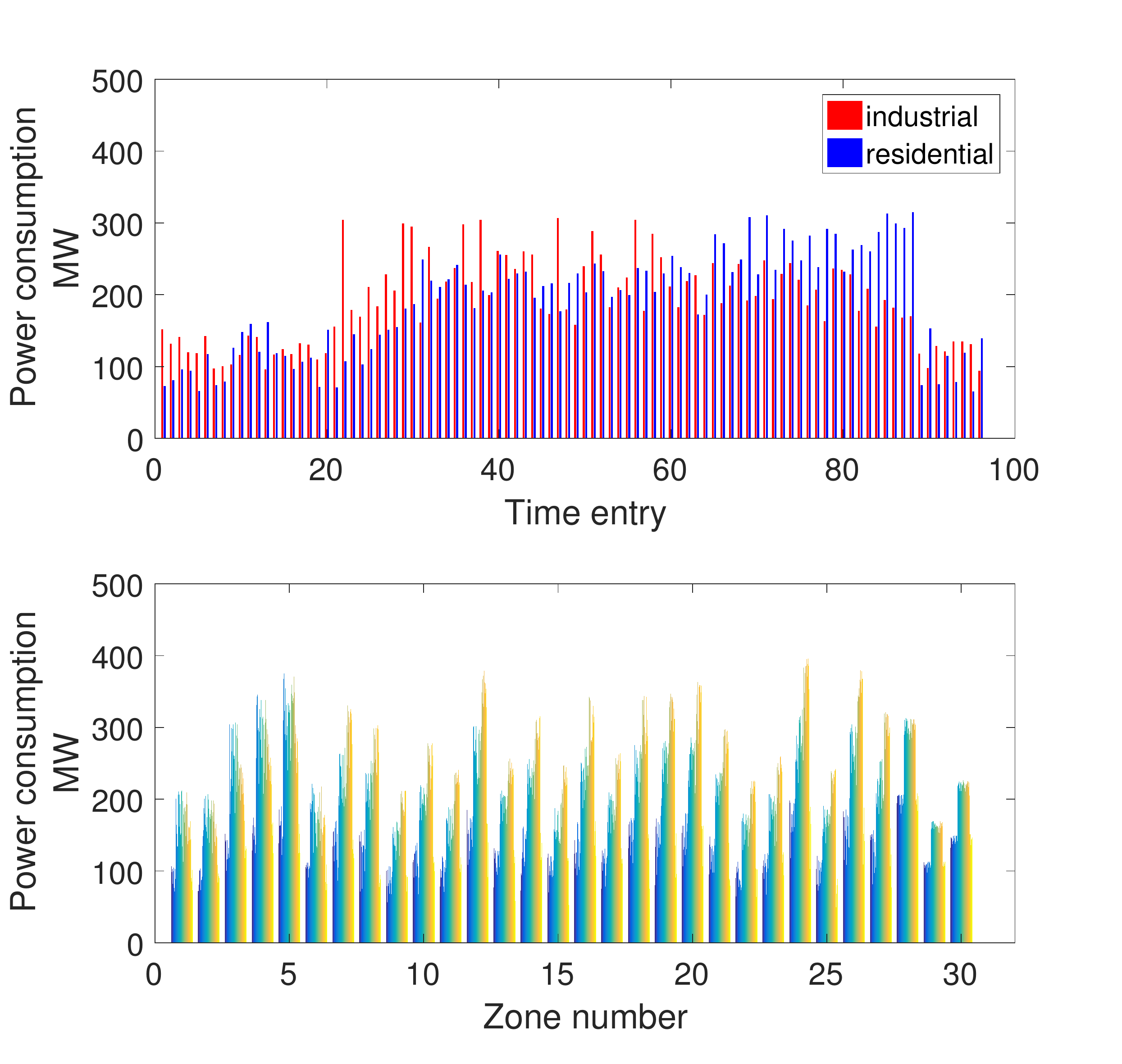}
\caption{Daily load profiles.}
\label{daily}
\vspace*{-0.4cm}
\end{figure}

The second part of our numerical tests involves making the actual schedules for all the zones. In this case, the number of outages and the duration are known. But electricity end users would prefer to know the exact time when load shedding is expected to happen so that they can make any necessary preparations. To figure out the information, the utility needs to have the load profiles for all the zones during the given period (30 days). 
All load profiles are generated to respect their distribution characteristics in practice. For example, high power consumption is expected for residential zones from 4pm to 8 pm as reflected in the profiles. The top plot of Figure \ref{daily} shows the daily load profile for one of the residential and industrial zones, namely zone 14 and zone 3, respectively. These two zones are selected because they have a similar average power $P^{\mathrm{ave}}_n$ but clearly different peak hours. The bottom plot shows the load profile for all the $N$ zones. Each profile is different because of the uniquely assigned value of $P^{\mathrm{ave}}_n$ and the randomness in distribution. Since the base time is 15 minutes we have 96 entries of data per zone for one day.

\begin{figure}[!t]
\centering
\includegraphics[width=3.4in]{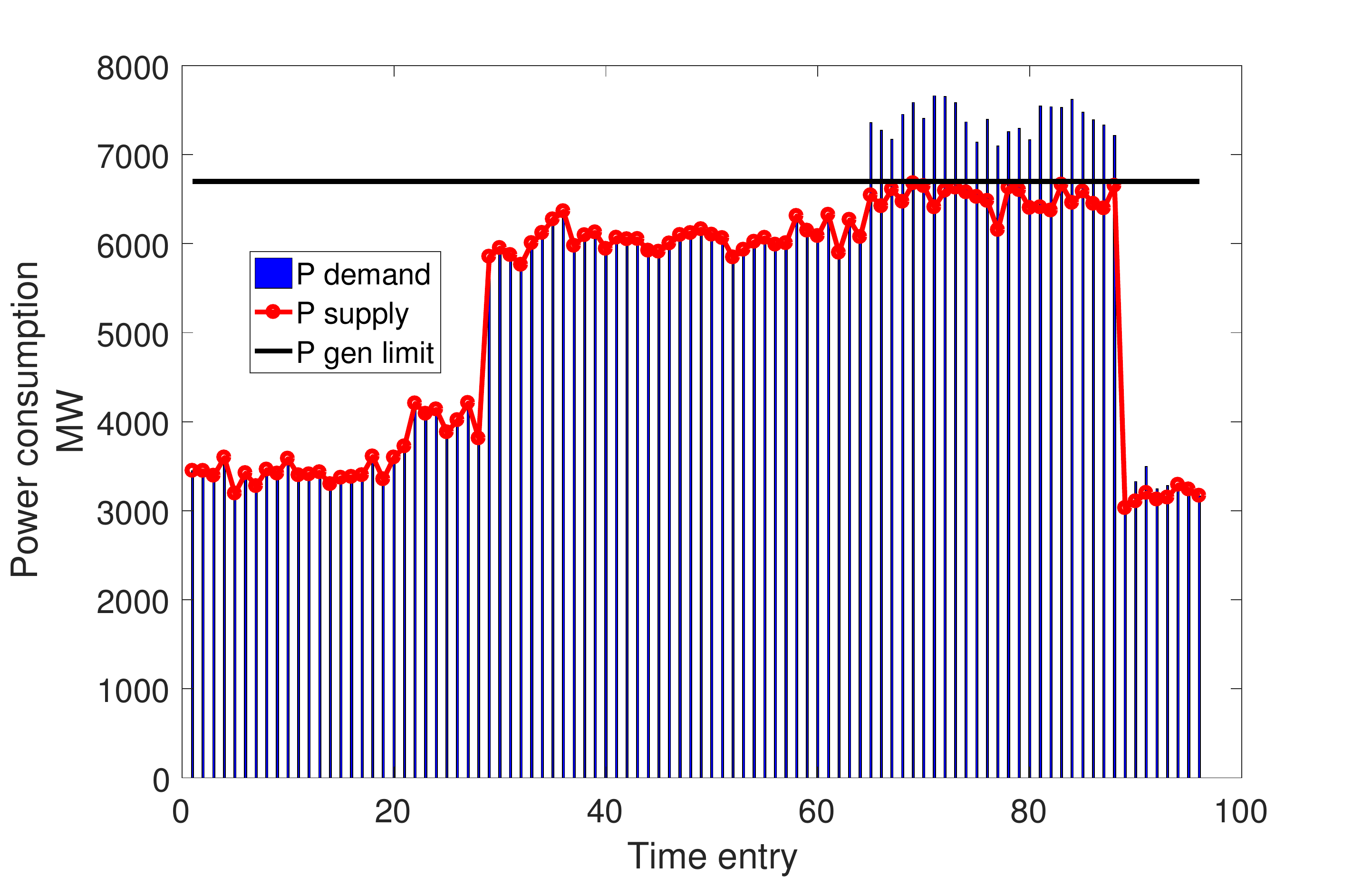}
\caption{Simulation results for load shedding over one day.}
\label{oneday}
\vspace*{-0.4cm}
\end{figure}

The max power generation capacity is selected to satisfy the $E_{\mathsf{sf}}$ specified in Table \ref{para}. In this case the power generation limit turns out to be $6700$ MW. If the system has solar or wind power generation, this maximum limit can be a time-varying parameter. However, for simplicity it is set to be a constant in our simulation. Figure \ref{oneday} shows the total power demand, where whenever it exceeds the power generation limit some load must be shed. The actual power being supplied is therefore always less than the generation limit.

\begin{figure}[!t]
\centering
\includegraphics[width=3.4in]{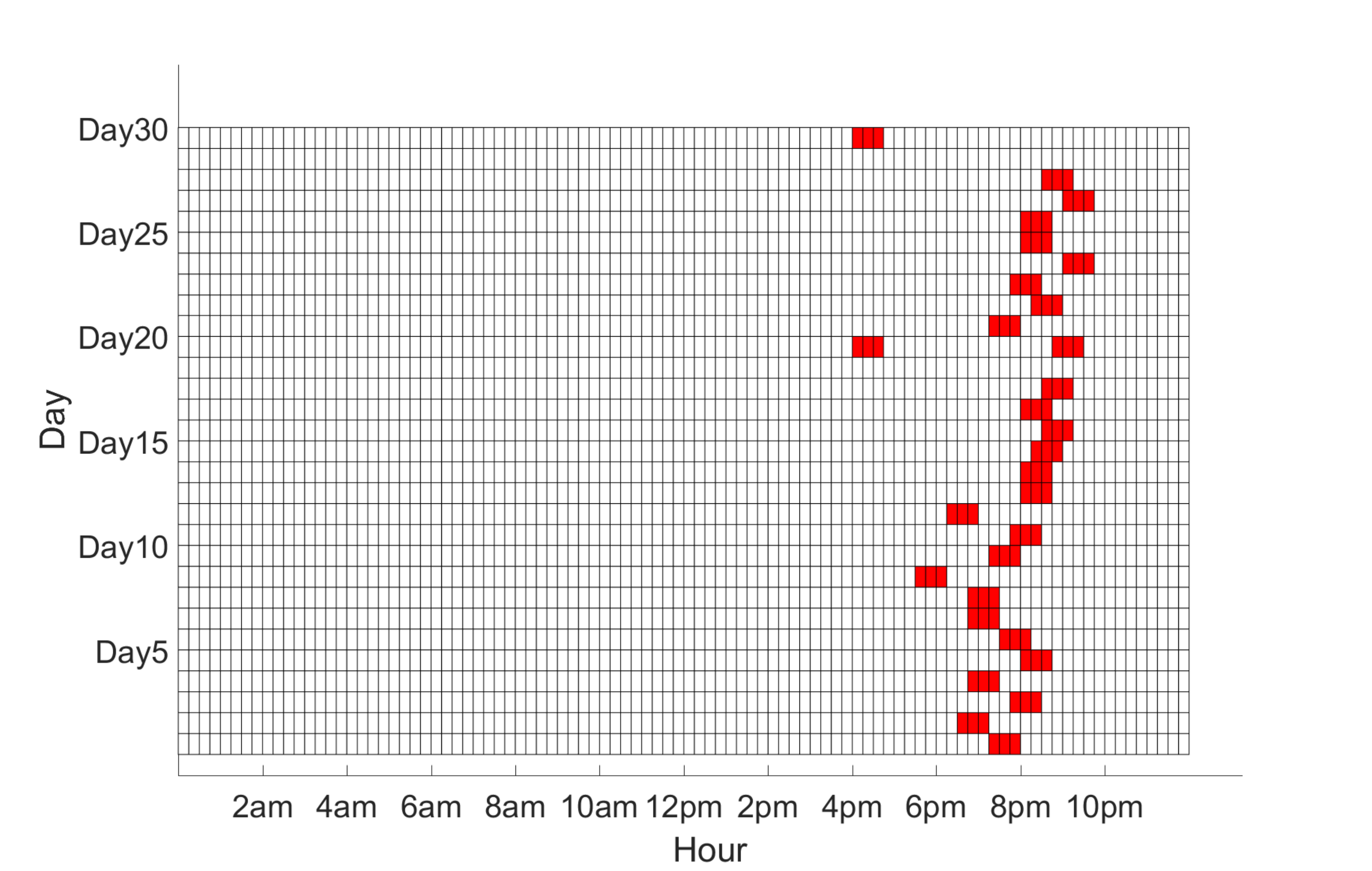}
\caption{Load shedding calendar for zone 10 for one month.}
\label{cal}
\vspace*{-0.4cm}
\end{figure}

Whenever the power demand exceeds the generation limit, the proposed model determines which load to shed based on the obtained optimal solution. This essentially creates a calendar of load shedding for each zone. As an example, consider the calendar shown in Figure \ref{cal} for the residential zone 10. From Figure \ref{result} it can be seen that this zone should have 29 outages in 30 days each of 45 minutes. Hence, a calendar can be created for each zone and provided to the respective consumers ahead of time.

\begin{figure}[!t]
\centering
\includegraphics[width=3.4in]{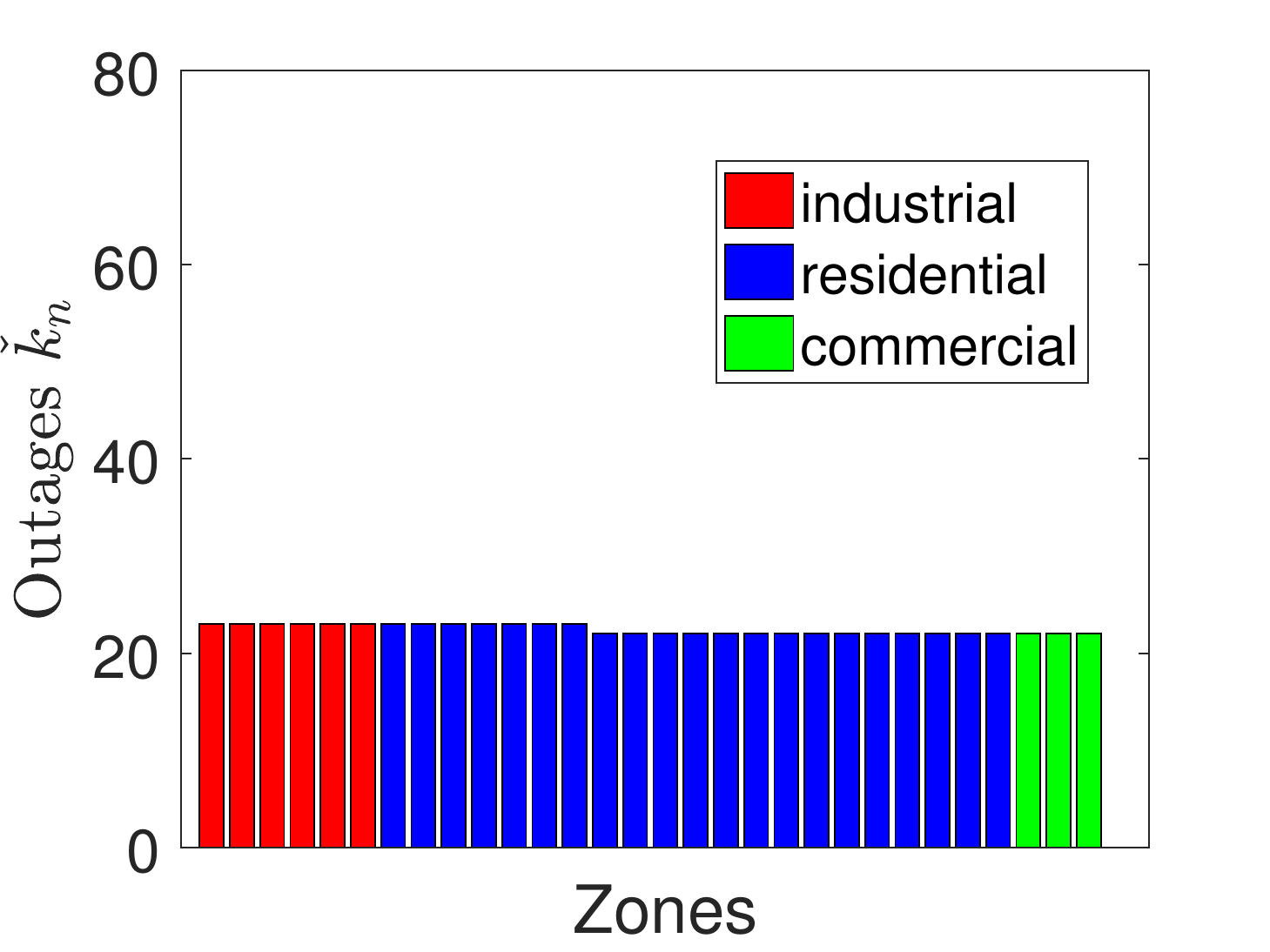}
\caption{Planned outages with the sequencing algorithm.}
\label{equal}
\end{figure}

\begin{figure}[!t]
\centering
\includegraphics[width=3.4in]{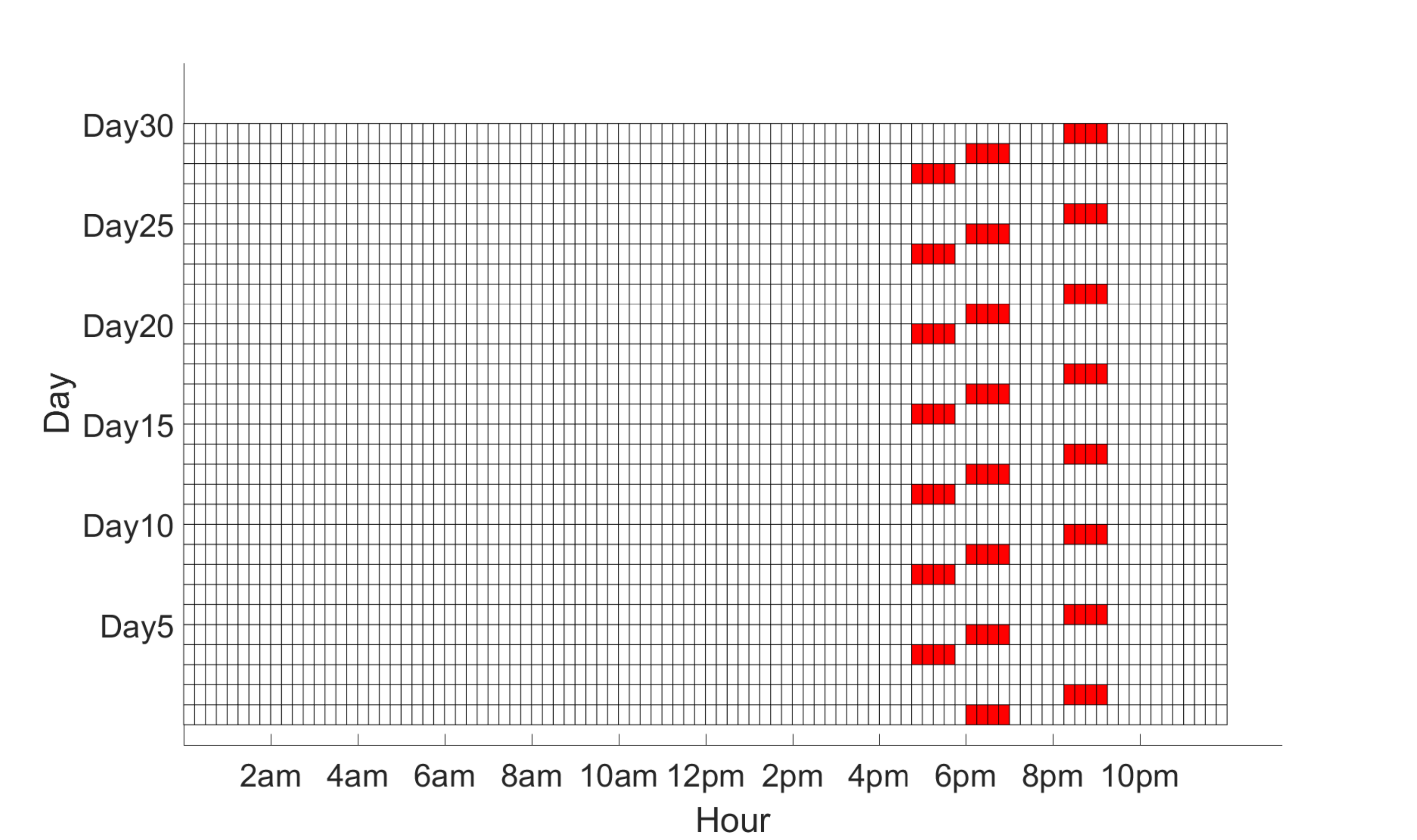}
\caption{Load shedding calendar for zone 10 with the sequencing algorithm.}
\label{cal2}
\vspace*{-0.4cm}
\end{figure}

Next, we will compare our proposed model and approach with some existing methods of load shedding. A common way to schedule rotational load shedding is based on sequencing the load shed across all the zones \cite{main}. This is the simplest way to ensure that all the zones equally share the outages. The \emph{round-robin sequence} $\{1,2,\ldots,N,1,2,...\}$ is chosen here while the duration for each outage is set to one hour. In other words, the first time when load has to be shed, zone 1 is shed for an hour. Then for zone 2, and so on and so forth. The resulting numbers of outages, and corresponding calendar for zone 10 are shown in Figure \ref{equal} and Figure \ref{cal2}, respectively. As expected the outages are distributed equally.  This calendar incidentally looks similar to the ones in Figure \ref{lesco1} and \ref{lesco2}. This fact implies that Lesco may implement the sequencing algorithm to distribute the outages among groups within a zone. With exactly the same setup (all parameters and load profiles), the load shedding cost of the sequencing algorithm turns out to be 465,640 units while the optimal cost of our proposed approach is 360,210 units, which yields a $29\%$ cost reduction.


\begin{figure}[!t]
\centering
\includegraphics[width=3.4in]{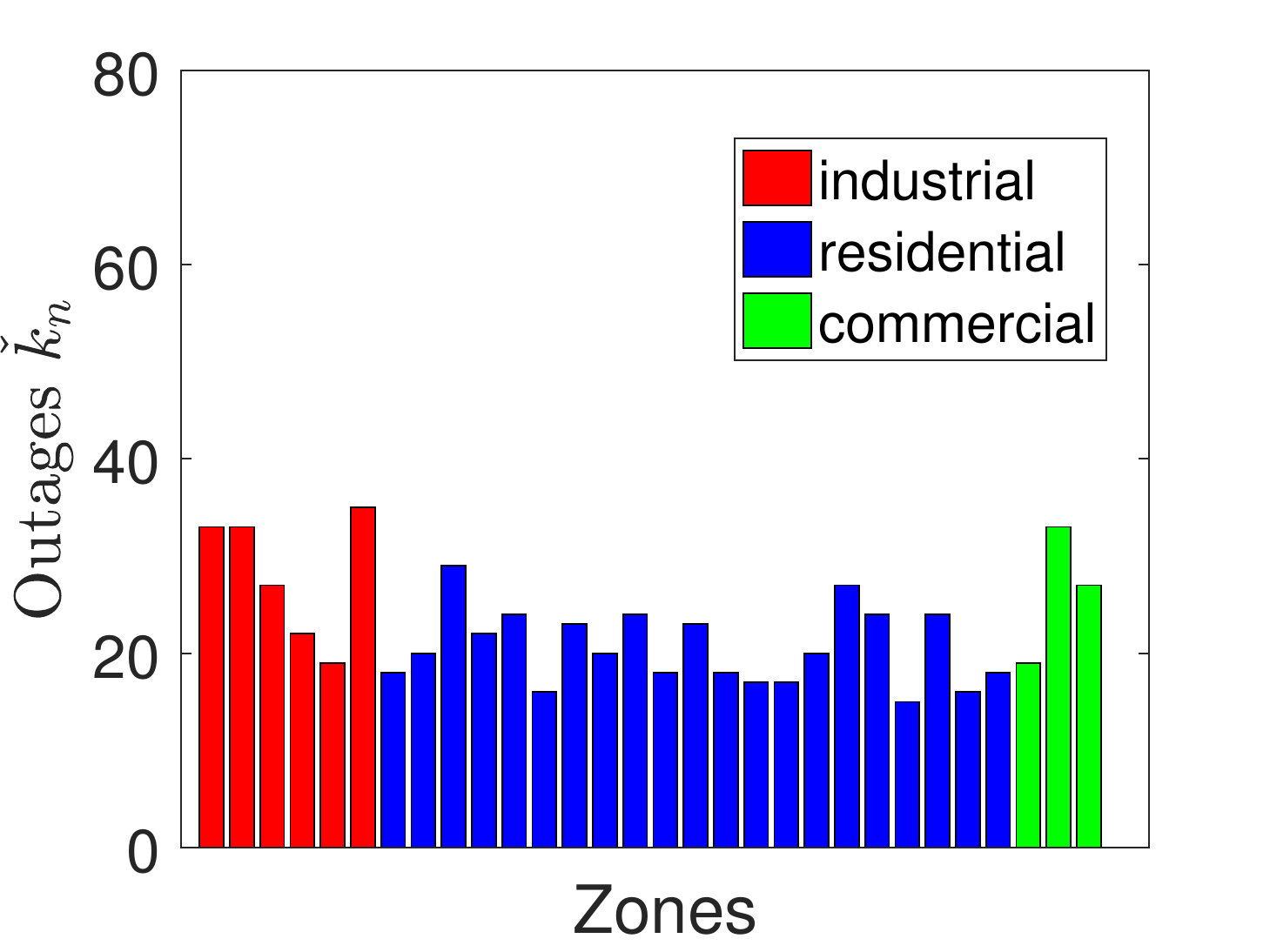}
\caption{Planned outages with the equal power shedding.}
\label{power}
\end{figure}

\begin{figure}[!t]
\centering
\includegraphics[width=3.4in]{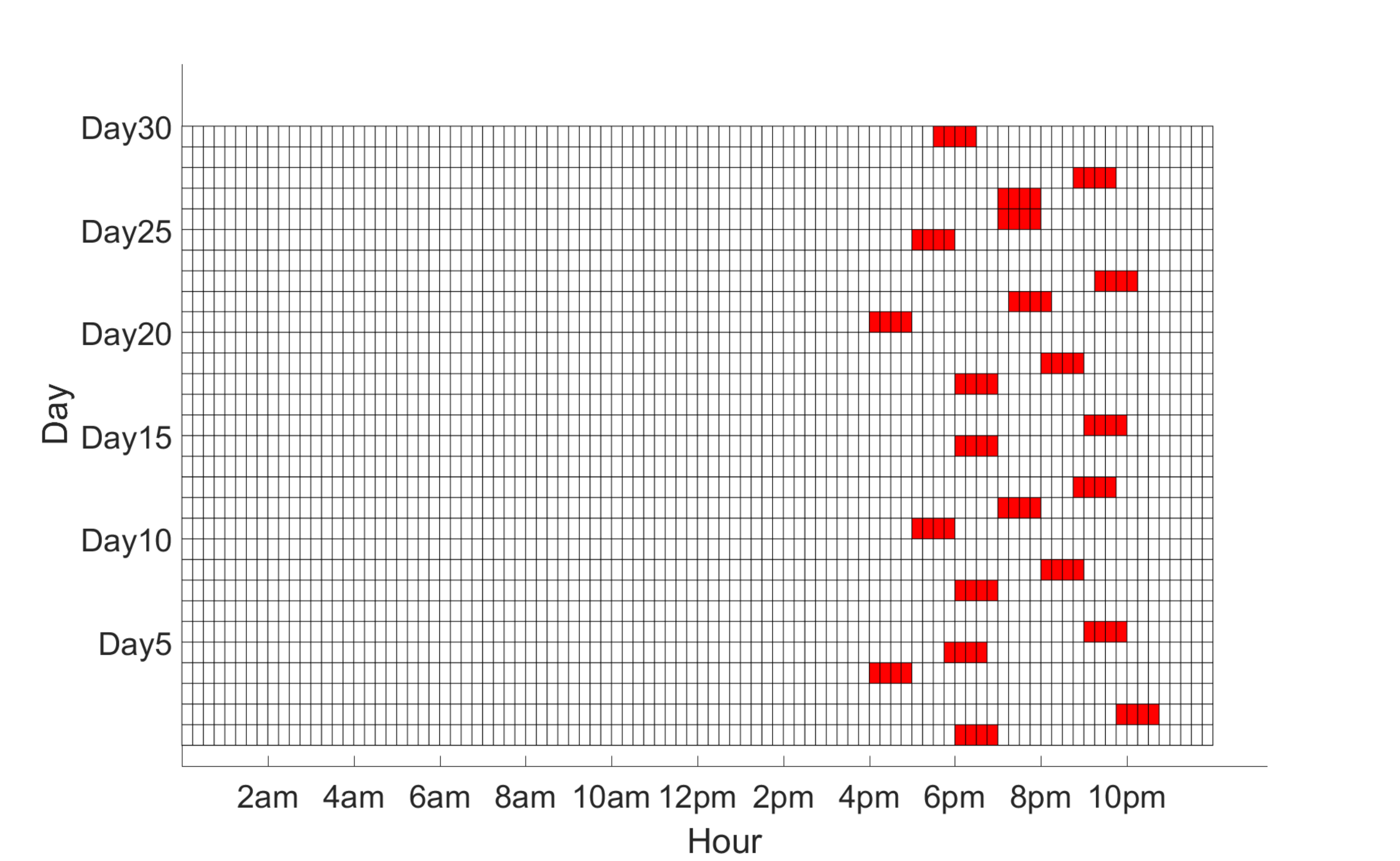}
\caption{Load shedding calendar for zone 10 with the equal power shedding.}
\label{cal3}
\vspace*{-0.1cm}
\end{figure}

Another  way to distribute power outages is the so-called \emph{equal power shedding}, which assigns outages such that each zone sheds roughly the same amount of power. Consequently, the zones with a larger power demand will have fewer planned outages. We simulate this method based on our setup, and it turns out that  each zone sheds about $21,086$ MWh of energy. The duration of each outage for all zones is still fixed to be one hour. Figure \ref{power} shows the resulting outages distribution while the load shedding calendar for zone 10 is given in Figure \ref{cal3}. Our proposed approach has $36\%$ cost reduction compared with this heuristic method.


\section{Conclusions}\label{sec:summary}
This paper explores optimal rotational load shedding in the presence of unavoidable power shortfalls. This is a challenging problem affecting millions of people every day, especially in a few developing countries. Variations of cost functions are proposed to model the heterogeneous preferences of load shedding for different zones. The suboptimal load shedding schedules are obtained by solving the formulated bilinear integer program via the McCormick relaxations. The numerical tests show the effectiveness and performance gain of the proposed approaches. Future research directions aim at handling the uncertainty of the estimates of the shortfall and average power consumption.

\balance

\bibliographystyle{IEEEtran}
\bibliography{loadshedding_bib}

\begin{thebibliography}{10}
\providecommand{\url}[1]{#1}
\csname url@samestyle\endcsname
\providecommand{\newblock}{\relax}
\providecommand{\bibinfo}[2]{#2}
\providecommand{\BIBentrySTDinterwordspacing}{\spaceskip=0pt\relax}
\providecommand{\BIBentryALTinterwordstretchfactor}{4}
\providecommand{\BIBentryALTinterwordspacing}{\spaceskip=\fontdimen2\font plus
\BIBentryALTinterwordstretchfactor\fontdimen3\font minus
  \fontdimen4\font\relax}
\providecommand{\BIBforeignlanguage}[2]{{%
\expandafter\ifx\csname l@#1\endcsname\relax
\typeout{** WARNING: IEEEtran.bst: No hyphenation pattern has been}%
\typeout{** loaded for the language `#1'. Using the pattern for}%
\typeout{** the default language instead.}%
\else
\language=\csname l@#1\endcsname
\fi
#2}}
\providecommand{\BIBdecl}{\relax}
\BIBdecl

\bibitem{power}
\BIBentryALTinterwordspacing
{Enerdata}, ``Global energy statistical yearbook 2017,'' 2018. [Online].
  Available: \url{https://yearbook.enerdata.net/}
\BIBentrySTDinterwordspacing

\bibitem{popu}
J.~Holdren, \emph{Population and environment}.\hskip 1em plus 0.5em minus
  0.4em\relax Kluwer Academic Publishers-Human Sciences Press, March 1991,
  vol.~12, ch. Population and the energy problem, pp. 231--255.

\bibitem{bitcoin}
\BIBentryALTinterwordspacing
{POWERCOMPARE}, ``Bitcoin mining now consuming more electricity than 159
  countries including {Ireland} \& most countries in {Africa},'' 2017.
  [Online]. Available: \url{https://powercompare.co.uk/bitcoin/}
\BIBentrySTDinterwordspacing

\bibitem{pak}
\BIBentryALTinterwordspacing
P.~Constable, ``A disaster in the making: {Pakistan’s} population surges to
  207.7 million,'' \emph{The Washington Post}, September 2017. [Online].
  Available: \url{https://wapo.st/2IQLmu9}
\BIBentrySTDinterwordspacing

\bibitem{pakfall}
\BIBentryALTinterwordspacing
K.~Kiani, ``Power cuts return as shortfall touches 7,000 mw,'' \emph{Dawn
  News}, May 2017. [Online]. Available: \url{https://www.dawn.com/news/1331738}
\BIBentrySTDinterwordspacing

\bibitem{impact}
S.~M. Hali, S.~Iqbal, W.~Yong, and S.~M. Kamran, ``Impact of energy sources and
  the electricity crisis on the economic growth: Policy implications for
  pakistan,'' \emph{Journal of Energy Tech. and Policy}, vol.~7, no.~2, 2017.

\bibitem{USIP}
R.~Aziz and M.~B. Ahmad, ``Pakistan’s power crisis,'' \emph{Special report.
  United States Institute of Peace. http://www. usip.
  org/sites/default/files/SR375-Pakistans-Power-Crisis-The-Way-Forward. pdf},
  2015.

\bibitem{lesco}
\BIBentryALTinterwordspacing
LESCO, ``Load shedding / shutdown schedule,'' 2018. [Online]. Available:
  \url{http://www.lesco.gov.pk/ShutdownSchedule}
\BIBentrySTDinterwordspacing

\bibitem{collapse}
T.~Alzahawi, M.~S. Sachdev, and G.~Ramakrishna, ``Time to voltage collapse
  dependence of load shedding to avoid voltage collapse,'' in \emph{Electrical
  Power and Energy Conf. (EPEC), 2017 IEEE}.\hskip 1em plus 0.5em minus
  0.4em\relax IEEE, 2017, pp. 1--5.

\bibitem{lesconews}
\BIBentryALTinterwordspacing
``Lesco announces 12-hour load-shedding schedule,'' \emph{Dunya News}, November
  2017. [Online]. Available:
  \url{http://dunyanews.tv/en/Pakistan/412805-LESCO-announces-12-hour-load-shedding-schedule}
\BIBentrySTDinterwordspacing

\bibitem{egypt}
\BIBentryALTinterwordspacing
K.~Fahim and M.~Thomas, ``Blackouts in {Egypt} prompt accusations,'' \emph{The
  New York Times}, August 2014. [Online]. Available:
  \url{https://www.nytimes.com/2014/08/23/world/middleeast/blackouts-in-egypt-prompt-accusations.html}
\BIBentrySTDinterwordspacing

\bibitem{ghana}
\BIBentryALTinterwordspacing
E.~Kwame, ``24-hr light-out in new loadshedding timetable?'' \emph{TV3Network},
  January 2015. [Online]. Available:
  \url{http://www.tv3network.com/business/24-hr-light-out-in-new-loadshedding-timetable}
\BIBentrySTDinterwordspacing

\bibitem{rsa}
\BIBentryALTinterwordspacing
C.~McGreal, ``Gold mines shut as south africa forced to ration power supply,''
  \emph{The Guardian}, January 2008. [Online]. Available:
  \url{https://www.theguardian.com/world/2008/jan/26/southafrica.international}
\BIBentrySTDinterwordspacing

\bibitem{indiafall}
A.~K. Arya, S.~Chanana, and A.~Kumar, ``Role of smart grid to power system
  planning and operation in india,'' in \emph{Proceedings of Intl. Conf. on
  Emerging Trends in Engr. and Tech.}, 2013, pp. 793--802.

\bibitem{chennai}
T.~Kogo, S.~Nakamura, S.~Pravinraj, and B.~Arumugam, ``A demand side prediction
  method for persistent scheduled power-cuts in developing countries,'' in
  \emph{IEEE PES Innovative Smart Grid Tech. Conf., Europe (ISGT-Europe)},
  2014, pp. 1--6.

\bibitem{pge}
\BIBentryALTinterwordspacing
{Pacific Gas and Electric Company}, ``Find your rotating outage block,'' 2018.
  [Online]. Available: \url{https://goo.gl/higJKb}
\BIBentrySTDinterwordspacing

\bibitem{japan}
\BIBentryALTinterwordspacing
P.~Ferguson, ``Tokyo region to face rolling blackouts through summer,''
  \emph{CNN}, March 2011. [Online]. Available:
  \url{www.cnn.com/2011/WORLD/asiapcf/03/31/japan.blackouts/index.html}
\BIBentrySTDinterwordspacing

\bibitem{dr}
S.~Khemakhem, M.~Rekik, and L.~Krichen, ``Optimal appliances scheduling for
  demand response strategy in smart home,'' in \emph{18th Intl. Conf. on Sci.
  and Tech. of Automatic Control and Computer Engr. (STA)}.\hskip 1em plus
  0.5em minus 0.4em\relax IEEE, 2017, pp. 546--550.

\bibitem{rsa2}
\BIBentryALTinterwordspacing
ESKOM, ``Interpreting {Eskom} load-shedding stages,'' 2018. [Online].
  Available:
  \url{http://loadshedding.eskom.co.za/LoadShedding/ScheduleInterpretation}
\BIBentrySTDinterwordspacing

\bibitem{pred}
P.~Chen, W.~Li, Y.~Chen, K.~Guo, and Y.~Niu, ``A parallel evolutionary extreme
  learning machine scheme for electrical load prediction,'' in \emph{Computing
  Conference, 2017}.\hskip 1em plus 0.5em minus 0.4em\relax IEEE, 2017, pp.
  332--339.

\bibitem{milp}
A.~Gupte, S.~Ahmed, M.~S. Cheon, and S.~Dey, ``Solving mixed integer bilinear
  problems using milp formulations,'' \emph{SIAM Journal on Optimization},
  vol.~23, no.~2, pp. 721--744, 2013.

\bibitem{main}
R.~Billinton and J.~Satish, ``Effect of rotational load shedding on overall
  power system adequacy indices,'' \emph{IEE Proceedings-Generation,
  Transmission and Distribution}, vol. 143, no.~2, pp. 181--187, 1996.

\end{thebibliography}
\end{document}